\documentclass[12pt]{amsart}

\newcommand\Ack{\medskip\noindent{\bf Acknowledgment.}\enspace}

\newcommand {\subplus}{\mathop{{\supset}\llap{\raise
0.5pt\hbox{\normalfont\small+}\hskip 0.5pt}}}

\newcommand {\supplus}{\mathop{{\subset}\llap{\raise
0.5pt\hbox{\normalfont\small+}\hskip 0.5pt}}}

\newcommand {\Cee}    {{\mathbb  C}}

\newcommand {\Zee}    {{\mathbb  Z}}

\newcommand {\fa}     {{\mathfrak{a}}}

\newcommand {\fas}    {{\mathfrak{as}}}

\newcommand {\fb}     {{\mathfrak{b}}}
\newcommand {\fc}    {{\mathfrak{c}}}

\newcommand {\fg}     {{\mathfrak{g}}}    %
\newcommand {\fgl}    {{\mathfrak{gl}}}  %
\newcommand {\fG}    {{\mathfrak{G}}}  %
\newcommand {\fh}     {{\mathfrak{h}}}

\newcommand {\fk}     {{\mathfrak{k}}}

\newcommand {\fle}    {{\mathfrak{le}}}
\newcommand {\fm}     {{\mathfrak{m}}}

\newcommand {\fo}     {{\mathfrak{o}}}

\newcommand {\fs}     {{\mathfrak{s}}}

\newcommand {\fsl}    {{\mathfrak{sl}}}
\newcommand {\fsle}   {{\mathfrak{sle}}}

\newcommand {\fsvect} {{\mathfrak{svect}}}

\newcommand {\fv}     {{\mathfrak{v}}}     %
\newcommand {\fvect}  {{\mathfrak{vect}}}   %

\newcommand {\cal} {\mathcal}

\newcommand {\cD}     {{\cal D}}

\def \opname#1#2%
  {\expandafter\newcommand \csname #1\endcsname {{\mathop{#2}\nolimits}}}

\newcommand{\rmname}[1]
  {\expandafter\newcommand \csname #1\endcsname {{\operatorname{#1}}}}

\newcommand{\rmnameii}[2]
  {\expandafter\newcommand \csname #1\endcsname {{\operatorname{#2}}}}

\rmname{act}
\rmname{Ad}
\rmname{Add}
\rmname{ad}
\rmname{Alt}
\rmname{alt}
\rmname{Ann}
\rmname{antidiag}
\rmname{Ber}
\rmname{ber}
\rmname{Br}
\rmname{card}
\rmname{ch}
\rmname{Char}
\rmname{cem}
\rmname{cj}
\rmname{Cliff}
\rmname{cntr}
\rmname{codim}
\rmname{coind}
\rmname{const}
\rmname{col}
\rmname{cork}
\rmname{cpr}
\rmname{diag}
\rmnameii{Div}{div}
\rmname{Def}
\rmname{Der}
\rmname{Dim}
\rmname{End}
\rmname{Even}
\rmname{Ext}
\rmname{gr}
\rmname{Hom}
\rmname{HT}
\rmnameii{Ht}{ht}
\rmname{hwt}
\rmname{Id}
\rmname{id}
\rmname{ind}
\rmname{Ind}
\rmname{Inf}
\rmname{irr}
\rmname{Le}
\rmname{Lie}
\rmname{lwt}
\rmname{mult}
\rmname{Mor}
\rmname{nm}
\rmname{Ob}
\rmname{Odd}
\rmname{Osc}
\rmname{per}
\rmname{Pic}
\rmname{pr}
\rmname{pro}
\rmname{Prime}
\rmname{Proj}
\rmname{prt}
\rmname{pt}
\rmname{Q}
\rmname{qet}
\rmname{qtr}
\rmname{rd}
\rmname{rk}
\rmname{row}
\rmname{Res}
\rmname{salt}
\rmname{Sch}
\rmname{SBr}
\rmname{scalar}
\rmname{Ser}
\rmname{sign}
\rmname{Smbl}
\rmname{spin}
\rmname{ssym}
\rmname{str}
\rmname{st}
\rmname{sgn}
\rmname{sq}
\rmname{symm}
\rmname{supp}
\rmname{Supp}
\rmname{St}
\rmname{Spec}
\rmname{Spm}
\rmname{tr}
\rmname{vpt}
\rmname{weyl}
\rmname{Weyl}
\rmname{Witt}

\opname{vvol}  {{v\hspace{-0.1ex}o\hspace{-0.02ex}l\/}}
\opname{pnt}  {\text{\normalfont pt}}
\opname{Span} {{Span}}
\opname{slim} {\overline{\lim}}
\opname{Vol}  {{V\hspace{-0.55ex}o\hspace{-0.02ex}l\/}}
\opname{QVol} {{Q\hspace{-0.3ex}V\hspace{-0.55ex}o\hspace{-0.02ex}l\/}}
\opname{PoVol}{{P\hspace{-0.35ex}o\hspace{-0.25ex}V\hspace{-0.55ex}o\hspace{-0.02ex}l\/}}
\opname{BVol} {{B\hspace{-0.2ex}V\hspace{-0.55ex}o\hspace{-0.02ex}l\/}}
\opname{Par}  {{P\hspace{-0.3ex}a\hspace{-0.05ex}r\/}}

\rmname{Mat}
\rmname{Bil}
\rmname{Diff}
\rmname{Ker}
\rmname{Herm}
\rmname{Coker}
\rmname{Conn}
\rmname{Covect}
\rmname{Vect}
\rmname{Int}

\rmnameii {IM} {Im}
\rmnameii {RE} {Re}

\opname{Aut} {{A\hspace{-0.2ex}u\hspace{-0.1ex}t\/}}
\opname{GL} {{G\hspace{-0.3ex}L}}
\opname{SL} {{S\hspace{-0.3ex}L}}
\opname{Exp} {{E\hspace{-0.2ex}x\hspace{-0.1ex}p\/}}
\opname{GQ} {{G\hspace{-0.2ex}Q}}
\opname{OSp} {{O\hspace{-0.25ex}S\hspace{-0.15ex}p\/}}
\opname{Out} {{O\hspace{-0.25ex}u\hspace{-0.15ex}t\/}}
\opname{Spp} {{S\hspace{-0.2ex}p\/}}
\opname{SpO} {{S\hspace{-0.2ex}p\hspace{-0.02ex}O\/}}
\opname{Pe} {{P\hspace{-0.25ex}e\/}}
\opname{SPe} {{S\hspace{-0.25ex}P\hspace{-0.25ex}e\/}}
\opname{Spin} {{S\hspace{-0.25ex}p\hspace{-0.05ex}i\hspace{-0.1ex}n\/}}
\opname{Iso} {{I\hspace{-0.25ex}s\hspace{-0.1ex}o\/}}
\opname{SSPe} {{S\hspace{-0.25ex}S\hspace{-0.15ex}P\hspace{-0.25ex}e\/}}
\opname{PeU} {{P\hspace{-0.25ex}e\hspace{-0.1ex}U\/}}
\opname{QU} {{Q\hspace{-0.15ex}U\/}}
\opname{U} {{U\/}}

\opname{cGQ} {{\cal G \hspace{-0.2em} Q \/}}
\opname{cSL} {{\cal S \hspace{-0.2em} L \/}}
\opname{cGL} {{\cal G \hspace{-0.2em} L \/}}
\opname{cOSp} {{\cal O \hspace{-0.2em} S \hspace{-0.3em} \it p\/}}
\opname{cPe} {{\cal P \hspace{-1.5pt} \it e\/}}
\opname{cVect} {{\cal V \hspace{-1.5pt} \it e\hspace{-0.1ex}c\hspace{-0.1ex}t\/}}
\opname{cVol} {{\cal V \hspace{-1.5pt} \it o\hspace{-0.1ex}l\/}}
\opname{cAut} {{\cal A \hspace{-0.2em} \it u\hspace{-0.1em}t\/}}
\opname{cCovect} {{\cal C \hspace{-1.5pt}
     \it o\hspace{-0.1ex}v\hspace{-0.1ex}e\hspace{-0.1ex}c\hspace{-0.1ex}t\/}}
\opname{CW} {{C\hspace{-0.15ex}W}}

\opname {Ab}   {{\sf Ab}}
\opname {Alg}   {{\sf Alg}}
\opname {ASch}  {{\sf Aff\;Sch}}
\opname {Funct}   {{\sf Funct}}
\opname {Gr}   {{\sf Gr}}
\opname {Grf}  {{{\sf Gr}_f}}
\opname {Mods}   {{\sf mods}}
\opname {Rings}   {{\sf Rings}}
\opname {Salg}   {{\sf Salg}}
\opname {Sets} {{\sf Sets}}
\opname {SSMan} {{\sf SMan}}
\opname {Top}  {{\sf Top}}
\opname {Vebun}   {{\sf Vebun}}

\newcommand {\ev} {{\bar0}}
\newcommand {\od} {{\bar1}}
\newcommand {\eps} {\varepsilon}

\newcommand {\pder}[1] {{\frac{\partial}{\partial {#1}}}}

\newcommand {\bcdot}   {\mathbin{\hbox{\raise.4ex\hbox{\bf.}}}} 

\newcommand {\secno} {}

\newtheorem{Theorem}{\secno Theorem}

\newenvironment {th*}[1]
    {\gdef\thname{#1} \begin{thn}}%
    {\end{thn}}
\newtheorem{thn}[Theorem] {\thname}

\theoremstyle{definition}

\newenvironment {ex*}[1]
    {\gdef\thname{#1} \begin{exn}}%
    {\end{exn}}
\newtheorem{exn}[Theorem]{\thname}

\theoremstyle{remark}

\newenvironment {rem*}[1]
    {\gdef\thname{#1} \begin{remn}}%
    {\end{remn}}
\newtheorem{remn}[Theorem]{\thname}

\newcommand {\ssec}{\subsection*}

\begin{document}

\title[Riemann tensor for exceptional supermanifolds]{The
analogs of the Riemann tensor for exceptional structures on
supermanifolds} \footnote{Appeared in:  S.~K.~Lando,
O.~K.~Sheinman (eds.)  {\it Proc.  International conference \lq\lq
Fundamental Mathematics Today" (December 26--29, 2001) in honor of
the 10th Anniversary of the Independent University of Moscow},
IUM, MCCME 2003, 89--109.}

\author{Pavel Grozman, Dimitry Leites, Irina Shchepochkina${}^1$}

\address{Department of Mathematics, University of Stockholm,
Roslagsv.  101, Kr\"aftriket hus 6, SE-104 05 Stockholm, Sweden;
mleites@math.su.se; ${}^1$On leave of absence from the Independent
University of Moscow, Bolshoj Vlasievsky per, dom 11, RU-121 002
Moscow, Russia; irina@mccme.ru}



\begin{abstract} H.~Hertz called any manifold $M$ with a given nonintegrable
distribution {\it nonholonomic}.  Vershik and Gershkovich stated
and R.~Montgomery proved that the space of germs of any
nonholonomic distribution on $M$  with an open and dense orbit of
the diffeomorphism group is either (1) of codimension one or (2)
an Engel distribution.

No analog of this statement for supermanifolds is formulated yet,
we only have some examples: our list (an analog of \'E.~Cartan's
classification) of simple Lie superalgebras of vector fields with
polynomial coefficients and a particular (Weisfeiler) grading
contains 16 series similar to contact ones and 11 exceptional
algebras preserving nonholonomic structures.

Here we compute the cohomology corresponding to the analog of the
Riemann tensor for the {\it super}manifolds corresponding to the
15 exceptional simple vectorial Lie superalgebras, 11 of which are
nonholonomic.  The cohomology for analogs of the Riemann tensor
for the manifolds with an exceptional Engel manifolds are computed
in \cite{L0}.

\end{abstract}

\maketitle

\section*{Introduction}

\ssec{The main result} In this paper the ground field is $\Cee$.
Here, for each of the 15 exceptional simple infinite dimensional
vectorial Lie superalgebras $\fg=\mathop{\oplus}\limits_{i\geq
-d}\fg_i$ in their Weisfeiler grading, we have computed
$H^i(\fg_-; \fg)$ for $i=0, 1, 2$ and
$\fg_-=\mathop{\oplus}\limits_{i<0}\fg_i$.  These cohomology are
especially interesting for $i=2$. If $d=1$, then $H^2(\fg_{-};
\fg)$ can be interpreted as the space of values of the  curvature
tensor for the $G$-structure, where $G$ is a Lie supergroup whose
Lie superalgebra is $\fg_0$; if $d>1$ we interprete $H^2(\fg_-;
\fg)$ as the space of values of the recently introduced {\it
nonholonomic curvature} for the nonintegrable distribution.

To make the text of interest to a wider audience, we would like to
review the list of simple vectorial Lie superalgebras and some not
so known background, cf. \cite{D}, but had to delete this from
this text for the lack of space; for the same reason we omitted
the results of computations of $H^i(\fg_-; \fg)$ for $i=0, 1$
($H^1(\fg_-; \fg)$ also has an interpretation of interest: its
elements represent derivations of $\fg_-$ into $\fg$, cf.
\cite{CK1}.)

The results demonstrate one more range of applicability of
SuperLie package. It is designed for various computations with Lie
superalgebras, not only for computation of (co)homology; for other
results and comparison with hand-made calculations, see, e.g.,
\cite{GLS2}. SuperLie is {\it Mathematica}-based which facilitates
its usage but imposes in-build {\it Mathematica} restrictions. We
hope to draw attention to possibilities SuperLie (now installed at
MPIM, Bonn; LPT, Ecole Normale Superior (Paris); Department of
Mathematics, University of Stockholm) reveals to its user.

In particular, our results (as well as similar results of
Poletaeva \cite{P1} --- \cite{P3} (now under one roof as
\cite{P4}) performed by bare hands) vividly demonstrate that in
the absence of complete reducibility computer-aided study is
indispensable.

\ssec{0.1. A result of R.~ Montgomery, Vershik and Gershkovich.
Nonholonomic curvature} R.~ Montgomery \cite{Mo} proved the
following statement whose particular case was stated in
\cite{VG1}. Let $W^k_n$ be the space of germs of $k$-dimensional
distributions at $0\in \Cee^n$. (Both \cite{Mo} and \cite{VG1}
consider the real case but the result is the same.) The group
$\text{Diff}_n$ of germs of diffeomorphisms of $\Cee^n$ acts on
$W^k_n$ and it is interesting to find out the conditions for
existence of the frame (i.e., point-wise values of the vector
fields from a given distribution) that generates a finite
dimensional (hence, nilpotent) Lie algebra. The answer:
\medskip

{\sl For $1<k<n-1$ and $(k, n)\neq (2, 4)$, any parametrization of
any open subset of the space of generic orbits of
$\text{Diff}_n$-action on $W^k_n$ requires $\geq k(n-k)-n$
functions in $n$ indeterminates. The exceptional case $k=1$ is
trivial. So  $W^k_n$  has an open and dense $\text{Diff}_n$-orbit
if and only if either (1) $k=n-1$ (for $n$ odd, this is the
contact structure) or (2) $(k, n)= (2, 4)$ (and then the
distribution is an Engel one)}.
\medskip

In \cite{L0}, the notion of  nonholonomic curvature is introduced
and the above cases (1) and (2) are considered. It turns out that
the nonholonomic curvature vanishes if $k=n-1$ (this is a
reformulation of Darboux theorem on a canonical form of the
contact form for $n$ odd) whereas for $(k, n)= (2, 4)$ --- the
Engel distribution --- $\dim H^2(\fg_-; \fg)=2$.

Observe that the (infinite dimensional) algebra of symmetries of
the exceptional (case (1) or (2)) nonholonomic distribution is
simple only if the distribution is of codimension 1 and $n$ is
odd. Contemporary mathematicians are often more than necessary
fixed on simple Lie algebras and this is, perhaps, an explanation
why the Lie algebra preserving an Engel structure (cf. \cite{L0})
was neglected for a long time. Differential geometers, though
mildly interested in the cases where the total algebra of
symmetries is simple, are more interested in cases where it is of
finite dimension (simple or not), and therefore their interests
are orthogonal to ours as is seen from motivations and results
reviewed, e.g., in \cite{Y}.

We do not know any super version of the above result of Vershik
and Gershkovich \cite{VG1} but we classified simple Lie
superalgebras of vector fields (\cite{LSh0}, \cite{LSh1},
\cite{LSh3}, cf. \cite{Ka7}) and, we see that, unlike non-super
case, there are 16 series and 11 exceptional simple vectorial Lie
superalgebras that preserve nonholonomic distributions. The series
will be considered elsewhere, \cite{GLS3}.

\ssec{0.2.  A nonholonomic analog of the Riemann tensor} H.~Hertz
called any manifold with a nonintegrable distribution a {\it
nonholonomic} one.  Until 1989, there was no general definition of
the analog of the Riemann tensor for nonholonomic manifolds, cf.
lamentations in \cite{VG1} and \cite{WB}, though all the
ingredients had been discovered (\cite{T}, \cite{Y}). Vershik even
conjectured  \cite{V} that such a general definition does not
exist, though in particular cases of small dimension the
nonholonomic curvature tensor was computed. In particular, in
supergravity.

Recall that SUGRA($N$) is a supergravity theory (or equations
thereof) on an $N$-extended Minkowski superspace. Whatever
SUGRA($N$) and Minkowski superspace are (there are several
versions of the definition and, unless $N=1$, there is no
consensus among physicists which of the definitions is \lq\lq it",
Manin \cite{M} suggests still other --- \lq\lq exotic" ---
versions of Minkowski superspaces, and this list of {\it ad hoc}
superizations of Minkowski space will, clearly, be continued, see
e.g., \cite{GL3}), they are superizations of the gravity theory
(or Einstein-Hilbert's equations) on the Minkowski space. So the
problem whose existence Wess honestly acknowledged in his lectures
\cite{WB} \lq\lq We do not know how to write the super Riemann
tensor" (on $N>2$ extended Minkowski superspaces) sounds strange:
take any textbook on differential geometry (say, \cite{St}) and
superize definition of the Riemann tensor or, more generally,
structure functions --- analogs of the Riemannian tensor --- for
any $G$-structure according to Sign Rule. This was precisely what
A. Schwarz with his colleagues and students  suggested to do
\cite{S}, see also \cite{CAF}.

The results of such an approach, however, seem to coincide with
the equations physicists write from their mysterious physical
considerations only for $N=1$  (but actually do not even in this
case, cf. \cite{S},  \cite{CAF} with \cite{GL1}).

In  \cite{GL1}, \cite{GL2} it was indicated that the roots of the
problem Wess addressed lie not in the  prefix \lq\lq super" which
only causes some signs in the classical definitions  of the
structure functions. The point is that every of numerous versions
of Minkowski {\bf super}space is {\bf nonholonomic}, unlike
Minkowski space, and since there was no general recipe for
computing nonholonomic analog of the curvature tensor, to write
SUGRA($N$) equations or even determine what is $N$-extended
Minkowski superspace (which of the versions satisfies some natural
requirements) was a problem.

Here we will reproduce the definition  (\cite{L0}, \cite{LP2},
\cite{GL1}) of the Riemann tensor $R$ in terms of Lie algebra
cohomology rather than in terms of Spencer cohomology (cf.
\cite{G}, \cite{St}) and give its generalization to nonholonomic
case.

Namely, fix a point on the nonholonomic (super)manifold with a
nonintegrable distribution $\cD$. Let $\cD$ be given by a system
of Pfaff equations and let $\fG$ be the filtered Lie superalgebra
preserving this system of  equations, $\fm$ the associated graded
one. Set $\fm_-=\mathop{\oplus}\limits_{i<0}\fm_i$; by default we
let $\fg_0=\fm_0$, the Lie superalgebra of grading preserving
derivations of $\fg_-=\fm_-$. Clearly, $\fm_-$ is nilpotent, see
\cite{VG2}, \cite{Y}. Very often an additional structure on $\cD$
is given; if this is the case, we take for $\fg_0$ a subalgebra of
$\fm_0$ that preserves this additional structure. Let $\fg$ be the
defined below generalized Cartan prolong of the pair $(\fg_-,
\fg_0)$.  Then, by the same arguments as in \cite{St}, {\bf the
possible values of the nonholonomic Riemann tensor $R$ at the
point span the superspace $H^2(\fg_-; \fg)$.}

\ssec{0.3.  The projective connections and their nonholonomic
analogs} The projective connection on the $n$-dimensional manifold
is the one whose group of automorphisms is locally isomorphic to
$\fsl(n+1)=\fg_-\oplus\fgl(n)\oplus(\fg_-)^*$, cf. \cite{G}. The
corresponding structure functions are from $H^2(\fg_-;
\fsl(n+1))$.

Similarly, for any $\Zee$-graded Lie superalgebra $\fg$ of finite
depth, let $\fh\subset\fg$ be a subalgebra with the same
nonpositive part. Then the elements of $H^2(\fg_-; \fh)$ are
analogs of the projective structure functions, especially
resembling them if $\dim \fh<\infty$. Such cohomology is
considered in \cite{GLS3}.

\section*{\S 1. Description of simple vectorial Lie
superalgebras}

For the lack of space we deleted all the preliminaries. The reader
willing to see them is referred to \cite{Sh14}, \cite{Sh5} and
\cite{LSh3}. For a detailed background see a preprinted version at
www.mpim-bonn.mpg.de. Observe only that $\Pi$ is the shift of
parity functor on superspaces, $\Vol$ is the space of densities
with the generator $\vvol$ in a fixed coordinate system.

\ssec{1.1.  The exceptional vectorial Lie subsuperalgebras} Here
are the terms $\fg_{i}$ for $i\leq 0$ of 14 of the 15 exceptional
algebras, the last column gives $\dim \fg_{-}$.  Here $\Lambda(n)$
is the Grassmann superalgebra with $n$ generators; $\id$ is the
identity representation of the subalgebra of the matrix Lie
superalgebra $\fgl(V)$ in the superspace $V$, let $\Lambda(\id)$
be the exterior algebra of $\id$; $\Vol_0$ is the space of
densities with integral 0;  and $T^0_{0}(\vec
0)=\Vol_0/\text{const}$ is well-defined only as module over
$\fsvect$: \footnotesize
$$
\renewcommand{\arraystretch}{1.3}
\renewcommand{\arraystretch}{1.4}\begin{tabular}{|c|c|c|c|c|}
\hline
$\fg$&$\fg_{-2}$&$\fg_{-1}$&$\fg_0$&$\dim\fg_{-}$\cr
\hline
\hline
$\fv\fle(4|3)$&$-$&$\Pi(\Lambda(3)/\Cee 1)$&$\fc(\fvect(0|3))$&$4|3$\cr
\hline
$\fv\fle(4|3; 1)$&$\Cee\cdot 1$&$\id\otimes\Lambda (2)$&
$\fc(\fsl(2)\otimes\Lambda(2)\subplus
T^{1/2}(\fvect(0|2))$&$5|4$\cr
\hline
$\fv\fle(4|3; K)$&$\id_{\fsl(3)}$&$\id_{\fsl(3)}\otimes \id_{\fsl(2)}\otimes
1$&$\fsl(3)\oplus\fsl(2)\oplus\Cee z$&$3|6$\cr
\hline
\hline
$\fv\fas(4|4)$&$-$&$\spin$&$\fas$&$4|4$\cr
\hline
\hline
$\fk\fas$&$\Cee\cdot 1$&$\Pi(\id)$& $\fc\fo(6)$&$1|6$\cr
\hline
$\fk\fas(; 1)$
&$\Lambda(1)$&$\id_{\fsl(2)}\otimes\id_{\fgl(2)}\otimes\Lambda(1)$&
$(\fsl(2)\oplus\fgl(2)\otimes\Lambda(1))\subplus\fvect(0|1)$&$5|5$\cr
\hline
$\fk\fas(; 3\xi)$&$-$& $\Lambda(3)$&$\Lambda(3)\oplus\fsl(1|3)$&$4|4$\cr
\hline
$\fk\fas(; 3\eta)$&$-$&$\Vol_{0}(0|3)$&
$\fc(\fvect(0|3))$&$4|3$\cr
\hline
\hline
$\fm\fb(4|5)$&$\Pi(\Cee\cdot 1)$&$\Vol (0|3)$&$\fc(\fvect(0|3))$&$4|5$\cr
\hline
$\fm\fb(4|5; 1)$&$\Lambda(2)/\Cee \cdot 1$
&$\id_{\fsl(2)}\otimes\Lambda(2)$
&$\fc(\fsl(2)\otimes\Lambda(2)\subplus T^{1/2}(\fvect(0|2))$&$5|6$\cr
\hline
$\fm\fb(4|5; K)$&$\id_{\fsl(3)}$&$\Pi(\id_{\fsl(3)}\otimes
\id_{\fsl(2)}\otimes
\Cee)$&$\fsl(3)\oplus\fsl(2)\oplus\Cee z$&$3|8$\cr
\hline
\hline
$\fk\fsle(9|6)$&$\Cee\cdot 1$&$\Pi(T^0_{0}(\vec 0))$&$\fsvect(0|4)_{3,
4}$&$9|6$\cr
\hline
$\fk\fsle(9|6; 2)$&$\id_{\fsl(3|1)}$&$\id_{\fsl(2)}\otimes\Lambda(3)$&
$\left (\fsl(2)\otimes\Lambda(3)\right) \subplus \fsl(1|3)$&$11|9$\cr
\hline
$\fk\fsle(9|6; K)$&$\id$&$\Pi(\Lambda^2(\id))$&$\fsl(5)$&$5|10$\cr
\hline
\end{tabular}
$$

\normalsize

\noindent Observe that none of the simple W-graded vectorial Lie
superalgebras is of depth $>3$ and only two algebras are of depth 3:
one of the above, $\fm\fb(4|5; K)$, for which we have $\fm\fb(4|5;
K)_{-3}\cong \Pi(\id_{\fsl(2)})$, and another one, $\fk\fsle(9|6;
CK)=\fc\fk(9|11)$.

This $\fc\fk((9|11)$ is the 15-th exceptional simple vectorial Lie
superalgebra; its non-positive terms are as follows (we assume that
the $\fsl(2)$- and $\fsl(3)$-modules are purely even):
$$
\renewcommand{\arraystretch}{1.4}
\begin{array}{l}
\fc\fk((9|11)_{0}\simeq \left
(\fsl(2)\oplus\fsl(3)\otimes\Lambda(1)\right) \subplus \fvect(0|1);\\
\fc\fk((9|11)_{-1}\simeq \id_{\fsl(2)}\otimes\left
(\id_{\fsl(3)}\otimes\Lambda(1)\right);\\
\fc\fk((9|11)_{-2}\simeq \id_{\fsl(3)}^*\otimes\Lambda(1);\\
\fc\fk((9|11)_{-3}\simeq \Pi(\id_{\fsl(2)}\otimes\Cee). \end{array}
$$

\ssec{1.2.  A description of $\fg$ as $\fg_{\ev}$ and $\fg_{\od}$}
In \cite{CK2} the exceptional algebras are described as
$\fg=\fg_{\ev}\oplus\fg_{\od}$.  For several series such
description is of little value because each homogeneous component
$\fg_{\ev}$ and $\fg_{\od}$ has a complicated structure.  For the
exceptions (and for twisted polyvector fields) the situation is
totally different!  Apart from being beautiful, such description
is useful for the construction of simple Volichenko algebras, cf.
\cite{LS}.

Recall in this relation a theorem \cite{Gr} that completely
describes bilinear differential operators acting in tensor fields
and invariant under all changes of coordinates.  It turned out
that almost all of the first order operators determine a Lie
superalgebra on its domain.  Some of these superalgebras are
simple or close to simple. In the constructions below we use some
of these invariant operators.

\noindent\underline{$\fg=\fk\fs\fle(5|10)$}:\; \;
$\fg_\ev=\fsvect(5|0)\simeq d\Omega^{4}$, \; \;
$\fg_\od=\Pi(d\Omega^{1})$ with the natural $\fg_\ev$-action on
$\fg_\od$ (the Lie derivative) and the bracketing of odd elements
being twice their product. (We identify:
$$
\partial_i=\sign(ijklm)dx_jdx_kdx_ldx_l\; \;  \text{
for any permutation $(ijklm)$ of $(12345)$}.
$$

\noindent\underline{$\fg=\fv\fas(4|4)$}:\; \; $\fg_\ev=\fvect(4|0)$,
and $\fg_\od=\Omega^{1}\otimes \Vol^{-1/2}$ with the natural
$\fg_\ev$-action on $\fg_\od$ and the bracketing of odd elements being
$$
[\omega_{1}\otimes \vvol^{-1/2}, \omega_{2}\otimes \vvol^{-1/2}]=
(d\omega_{1}\wedge\omega_{2}+ \omega_{1}\wedge d\omega_{2})\otimes
\vvol^{-1}],
$$
where we identify
$$
dx_{i}\wedge dx_{j}\wedge dx_{k}\otimes
\vvol^{-1}=\sign(ijkl)\partial_l\text{ for any permutation $(ijkl)$
of $(1234)$}.
$$

\noindent\underline{$\fg=\fv\fle(3|6)$}:\; \;
$\fg_\ev=\fvect(3|0)\oplus \fsl(2)^{(1)}_{\geq 0}$, where
$\fg^{(1)}_{\geq 0}=\fg\otimes\Cee[t]$ , and with the natural
$\fg_\ev$-action on $\fg_\od=\left(\Omega^{1}\otimes
\Vol^{-1/2}\right )\otimes \id_{\fsl(2)^{(1)}_{\geq 0}}$.

Recall that $\id_{\fsl(2)}$ is the irreducible $\fsl(2)$-module
$L^1$ with highest weight 1; its tensor square splits into
$L^2\simeq \fsl(2)$ and the trivial module $L^0$; accordingly,
denote by $v_{1}\wedge v_{2}$ and $v_{1}\bullet v_{2}$ the
projections of $v_{1}\otimes v_{2}\in L^1\otimes L^1$ onto the
skew-symmetric and symmetric components, respectively.  For
$f_{1}, f_{2}\in\Omega ^0$, $\omega_{1}, \omega_{2}\in\Omega ^1$
and $v_{1}, v_{2}\in L^1$, we set
$$
\renewcommand{\arraystretch}{1.4}
\begin{array}{l}
{}[(\omega_{1}\otimes v_{1})\vvol^{-1/2}, (\omega_{2}\otimes
v_{2})\vvol^{-1/2}]=\\
\left(\omega_{1}\wedge \omega_{2})\otimes
(v_{1}\wedge v_{2})+d\omega_{1}\wedge \omega_{2}+\omega_{1}\wedge
d\omega_{2})\otimes (v_{1}\bullet v_{2})\right)\vvol^{-1},
\end{array}
$$
where we identify $\Omega ^0$ with $\Omega ^3\otimes_{\Omega ^0}
\Vol^{-1}$ and $\Omega ^2\otimes_{\Omega ^0}
\Vol^{-1}$ with $\fvect(3|0)$ by setting
$$
dx_{i}\wedge dx_{j}\otimes
\vvol^{-1}=\sign(ijk)\pder{x_{k}}\text{ for any permutation $(ijk)$
of $(123)$}.
$$

\noindent\underline{$\fg=\fm\fb(3|8)$}:\; \;
$\fg_\ev=\fvect(3|0)\oplus \fsl(2)^{(1)}_{\geq 0}$, and
$\fg_\od=\fg_{-1}\oplus \fg_1$, where
$$
\fg_{-1}=
\left(\Pi\Vol^{-1/2}\right)\otimes \id_{\fsl(2)^{(1)}_{\geq 0}}\; \text{ and
}\; \fg_1= \left(\Omega^{1}\otimes \Vol^{-1/2}\right)\otimes
\id_{\fsl(2)^{(1)}_{\geq 0}};
$$
clearly, one can interchange $\fg_{\pm 1}$.

Multiplication is similar to that of $\fg=\fv\fle(3|6)$.  For
$f_{1}, f_{2}\in\Omega ^0$, $\omega_{1}, \omega_{2}\in\Omega ^1$
and $v_{1}, v_{2}\in L^1$, we set
$$
\renewcommand{\arraystretch}{1.4}
\begin{array}{l}
{}[(\omega_{1}\otimes v_{1})\vvol^{-1/2}, (\omega_{2}\otimes
v_{2})\vvol^{-1/2}]= 0, \\
{}[(f_{1}\otimes v_{1})\vvol^{-1/2}, (f_{2}\otimes v_{2})\vvol^{-1/2}]=
(df_{1}\wedge df_{2})\otimes (v_{1}\wedge v_{2})\vvol^{-1}, \\
{}[(f_{1}\otimes v_{1})\vvol^{-1/2}, (\omega_{1}\otimes v_{2})\vvol^{-1/2}]=
\left(f_{1}\omega_{1}\otimes (v_{1}\wedge v_{2})+(df_{1}\omega_{1}+
f_{1}d\omega_{1})\otimes (v_{1}\bullet v_{2})\right)\vvol^{-1}.
\end{array}
$$

\noindent\underline{$\fg=\fk\fas$}:\; \; $\fg_\ev=\fvect(1|0)\oplus
\fsl(4)^{(1)}_{\geq 0}$, and $\fg_\od=\fg_{-1}\oplus \fg_1$, where
$\fg_{-1}= \Pi\left(\Lambda ^2(\id_{\fsl(2)^{(1)}_{\geq 0}})\right)$ and
$\fg_1= \Pi\left(S^{2}(\id_{\fsl(2)^{(1)}_{\geq 0}})\right)$; clearly, one
can interchange $\fg_{\pm 1}$.

\section*{\S 2. Main result}

The above description of the exceptional algebras is nice to
visualize them, but in calculations we have used sometimes the
description of the elements from \cite{Sh14}, which we do not
reproduce to save space.

\begin{Theorem} The $\fg_0$-modules  $H^2(\fg_-; \fg)$
are given by the lists of modules over the semisimple part of $(\fg_0)_\ev$ given in \S 3. \end{Theorem}

{\bf Comments}. 1) {\bf important}: observe that if the
$\fg_0$-module from  $H^2(\fg_-; \fg)$ is indecomposable, then, in
equations like Einstein equations or Wess-Zumino constraints,
there is no need to vanish all its irreducible components: it
suffices to vanish only the modules that generate all.

2) The highest weights are given with respect to the standard
basis of Cartan subalgebra of the maximal semisimple part of
$(\fg_0)_\ev$, or, if this part is $\fsl(n)$, with respect to the
matrix units $E_{ii}$ of the $\fgl(n)$. The degree of the cocycle
is given relative to the grading of $\fg$. The weight of vector
$A$ is denoted by $w(A)$.
 \lq\lq mult" denotes the multiplicity of the corresponding module in
the space of closed/exact forms. If it is not equal to $r/0$, the
respective cohomologies are not pure (are defined up to a
coboundary) and it may well happen that our choice of the
representative may be beautified. The expression $d[v]^2$ means
$(dv)^ 2$. Few more notations appear in respective places.

3) For $\fg$ of depth $d$, the degrees of cohomology (interpreted
as orders $k$ of the structure functions responsible for
obstructions to flattening the structure under investigation up to
$k$th infinitesimal neighborhood, cf. \cite{St}) may range from
$2-d$, which we will indicate for $d>1$. Recall that the structure
functions of order $k$ are defined provided structure functions of
lesser orders vanish. It happens sometimes (the case of Riemannian
metric) that in these lesser orders there are no cohomology
(torsion-free property of Levi-Civita connection). Otherwise (if
something in lesser orders is nonzero) the vanishing conditions on
these low-order structure functions are analogs of Wess-Zumino
constraints in supergravity, cf. \cite{WZ}.

4) Due to a theorem formulated only for $d=1$ (by Serre for Lie
algebras and by Serganova for superalgebras, cf. \cite{LPS}), the
order of nonzero structure functions from $H^2(\fg_{-1}; \fg)$ is
always equal to 1 under certain conditions ({\it involutivity}).
As we will see, the considered Lie superalgebras of depth 1  are
all involutive. The yet nonexisting analog of Serre's theorem for
$d>1$ is more complicated: cohomology may be non-vanishing in
several degrees. If the lowest degrees are absent, we indicate
this by writing \lq\lq torsion-free" by analogy with the
Riemannian manifolds.

\section*{\S 3. The $\fg_0$-modules  $H^2(\fg_-; \fg)$}

\ssec{\protect $\fg=\fv\fle (4|3)$}

Cohomology: a single irreducible $\fg_0$-module in $\deg =1$, 
$\dim=(24|24)$.

Set:
$$
\text{$w(u_i)=(1,1,1)-\eps_i$, $w(y)=(0,0,0)$, $w(\xi_i)=\eps_i$
and $\partial_0=\pder{y}$, $\partial_i=\pder{u_i}$,
$\delta_i=\pder{\xi_i}$.}
$$
The $\fgl(3)$-highest vectors are:

\vskip 0.2 cm

{\small
\renewcommand{\arraystretch}{1.4}\begin{tabular}{|c|c|c|c|}
\hline
$\fgl(3)$-highest vectors & $\fgl(3)$-weight & dim & mult \\
\hline
$\partial_1\, d[\partial_2] \wedge d[\partial_3]$ & $(2, 0, 0)$ & $(6|0)$ & 3/2 \\
\hline
$\partial_0\, d[\partial_0] \wedge d[\delta_{1}]$ & $(1, 0, 0)$ & $(0|3)$ & 5/4 \\
\hline
$\partial_1\, d[\partial_3] \wedge d[\delta_{1}]$ & $(2, 0, -1)$ & $(0|15)$ & 2/1 \\
\hline
$\partial_1\, d[\partial_3] \wedge d[\partial_0]$ & $(1, 0, -1)$ & $(8|0)$ & 3/2 \\
\hline
$\partial_1\, d[\delta_{1}] \wedge d[\delta_{1}]$ & $(2,-1,-1)$ & $(10|0)$ & 1/0 \\
\hline
$\partial_1\, d[\partial_0] \wedge d[\delta_{1}]$ & $(1,-1,-1)$ & $(0|6)$ & 1/0 \\
\hline
\end{tabular}}

\ssec{\protect $\fg=\fv\fle (4|3; 1)$} $2-d=0$.

Cohomology in $\deg =1$, torsion-free, $\dim=(20|20)$:

The $(\fg_0)_\ev$-highest weight vectors (with respect to $\fsl(2)
\oplus \fgl(2)$)  are as follows: (the first coordinate of the
weight is given  with respect to a copy of $\fsl(2)$  realized as
$\fo(3)$, with half-integer weights; the last two coordinates are
with respect to a copy of $\fgl(2)$). In this realization
$$
\renewcommand{\arraystretch}{1.4}\begin{array}{l}
\text{$w(u_1)=(1,1,1)$, $w(u_2)=(1,0,1)$, $w(u_3)=(1,1,0)$,
$w(y)=(-1,0,0)$,}\\
\text{$w(\xi_1)=(0,0,0)$, $w(\xi_2)=(0,1,0)$, $w(\xi_3)=(0,0,1)$.}
\end{array}
$$
Denote the elements of $\fg_-$  as follows:
$$
\renewcommand{\arraystretch}{1.4}\begin{array}{l}
\text{${g_1}={\partial_2}$; ${g_2}={\partial_3}$;
${g_3}={\delta_2}$; ${g_4}={\delta_3}$;}\\
\text{ ${g_5}=-{u_2}{\partial_1}+{\xi_1}{\delta_2}$;
${g_6}=-{u_3}{\partial_1}+{\xi_1}{\delta_3}$;}\\
\text{
${g_7}=-{y}{\delta_3}-{\xi_1}{\partial_2}+{\xi_2}{\partial_1}$;
${g_8}={y}{\delta_2}-{\xi_1}{\partial_3}+{\xi_3}{\partial_1}$;
${g_9}={\partial_1}$.}
\end{array}
$$ We have

\vskip 0.2 cm

{\small

\renewcommand{\arraystretch}{1.4}\begin{tabular}{|c|c|c|c|c|}
\hline
N & $\fsl(2) \oplus \fgl(2)$-highest vectors & $\fsl(2) \oplus \fgl(2)$-weight & dim & mult \\
\hline
$[1]$ & $2\, {\partial_2}{{dg}_1}\wedge {{dg}_5}+{\partial_2}{{dg}_2}\wedge {{dg}_6}+
  {\partial_3}{{dg}_2}\wedge {{dg}_5}$ & $(0, 1, 0)$ & $(2|0)$ & 5/4 \\
\hline
$[2]$ & $\partial_2 {{dg}_2}\wedge{{dg}_5}-$ & $(0, 2, -1)$ & $(4|0)$ & 2/1 \\
\hline
$[3]$ & $-{\delta_3}{{dg}_3}\wedge {{dg}_8}+ ({y}{\delta_3} -
{\xi_1}{\partial_2}+{\xi_2}{\partial_1}){{dg}_8}\wedge{{dg}_8}$ & $(1,2,-1)$ & $(0|8)$ & 3/2 \\
\hline
$[4]$ & $\partial_2 {{dg}_5}\wedge{{dg}_6}$ & $(0, 1, 0)$ & $(2|0)$ & 1/0 \\
\hline
$[5]$ & ${2\partial_2}{{dg}_5}\wedge {{dg}_7}+{\partial_2}{{dg}_6}\wedge {{dg}_8}+
  {\partial_3}{{dg}_5}\wedge {{dg}_8}$ & $(1, 1, 0)$ & $(0|4)$ & 2/1 \\
\hline
$[6]$ & ${\partial_2}{{dg}_5}\wedge {{dg}_8}$ & $(1, 2, -1)$ & $(0|8)$ & 1/0 \\
\hline
$[7]$ & $-{\delta_3}{{dg}_5}\wedge {{dg}_8}-{\partial_2}{{dg}_8}\wedge {{dg}_8}$ & $(2,2,-1)$ & $(12|0)$ & 1/0 \\
\hline
\end{tabular}}

\vskip 0.2 cm

$\fg_0$-modules:

\vskip 0.2 cm

$[A]=[2]+[3]+[6]+[7]$ of $\dim=(16|16)$ and
$[B]=[A]+[1]+[4]+[5]$ of $\dim=(20|20)$; the modules $[A]$ and $[B]/[A]$ are irreducible.

\ssec{\protect $\fg=\fv\fle (4|3; K)$} $2-d=0$.

The $\fg_0$-highest weight vectors are as follows:
$\fg_0=\fsl(2) \oplus \fgl(3)$; the first coordinate of the weight is given  with respect to $\fsl(2)$  realized as
$$
\text{$x_-=\partial_y$, \quad $x_+=y^2\partial_y +
y\sum_i\xi_i\delta_i +\xi_1\xi_2\partial_3 - \xi_1\xi_3\partial_2
+ \xi_2\xi_3\partial_1$;}
$$
$\fgl(3)$ is realized as
$$
\text{$x^i_j=-u_j\partial_i+\xi_i\delta_j$ for $i\neq j$ and
$x^i_i=-u_i\partial_i+\xi_i\delta_i+\sum_k u_i\partial_i$.}
$$

In this realization
$$
\renewcommand{\arraystretch}{1.4}
\begin{array}{l}
\text{$w(u_1)=(2,0,1,1)$, $w(u_2)=(2,1,0,1)$,
$w(u_3)=(2,1,1,0)$,}\\
\text{$w(y)=(-2,0,0,0)$, $w(\xi_1)=(0,-1,0,0)$,
$w(\xi_2)=(0,0,-1,0)$, $w(\xi_3)=(0,0,0,-1)$).}
\end{array}
$$
\vskip 0.2 cm

Cohomology: a single irreducible $\fg_0$-module in $\deg=0$, $\dim=(30|0)$:

\nopagebreak
\vskip 0.2 cm

{\small

\renewcommand{\arraystretch}{1.4}\begin{tabular}{|c|c|c|c|}
\hline
$\fsl(2) \oplus \fgl(3)$-highest vectors & $\fsl(2) \oplus \fgl(3)$-weight & dim & mult \\
\hline
${\partial_1}{d[\delta_1]}^2$ & $(2,2,-1,-1)$ & $(30|0)$ & 1/0 \\
\hline
\end{tabular}}

\vskip 0.2 cm

\ssec{\protect $\fg=\fv\fa\fs(4|4)$}

The weights are given relative to $\fgl(4)\subset\fg_0$. In this
realization, the weight of $u_i$ is $\eps_i-\frac12(1,1,1,1)$ and
the weight of $\xi_i$ is $-\eps_i$; we set
$\partial_i=\pder{u_i}$, $\delta_i=\pder{\xi_i}$ for $1\le i\le
4$.

\vskip 0.2 cm

\nopagebreak

Cohomology: in $\deg= 1$, $\dim=(40|40)$:

\nopagebreak

\vskip 0.2 cm

{\small
\renewcommand{\arraystretch}{1.4}\begin{tabular}{|c|c|c|c|c|}
\hline
N& $\fgl(4)$-highest vectors & weight & dim & mult \\
\hline
$[1]$ & ${\delta_1}d[\partial_2]\wedge d[\partial_3]-{\delta_2}d[\partial_1]\wedge d[\partial_3]+
{\delta_3}d[\partial_1]\wedge d[\partial_2]$ & $(0, 0, 0, -1)$ & $(0|4)$ & 4/3 \\
\hline
$[2]$ & ${\partial_4}d[\partial_1]\wedge d[\partial_2]$ & $(\frac12, \frac12, -\frac12, -\frac32)$ & $(20|0)$ & 3/2 \\
\hline
$[3]$ & ${\partial_4}d[\delta_4]\wedge d[\delta_4]$ & $(\frac12, \frac12, \frac11, -\frac52)$ & $(20|0)$ & 1/0 \\
\hline
$[4]$ & ${\delta_1}d[\delta_4]\wedge d[\delta_4]$ & $(1, 0, 0, -2)$ & $(0|36)$ & 2/1 \\
\hline
\end{tabular}}

\vskip 0.2 cm

The module is reducible but indecomposable.  Vector $[1]$ is the
highest weight vector in the quotient module.

\ssec{\protect$\fg=\fk\fa\fs$} $2-d=0$. For brevity, instead of $K_f$ we write simply $f$ and $df^2$ means $(df)^2$.

\nopagebreak

Cohomology: a single irreducible module in $\deg= 1$, hence torsion free. The $\fo(6)$-weights are:

\vskip 0.2 cm

\nopagebreak

{\small
\renewcommand{\arraystretch}{1.4}\begin{tabular}{|c|c|c|c|}
\hline
$\fo(4)$-highest vectors & $\fo(6)$-weight & dim & mult \\
\hline
$\xi_1 \xi_2  d[\xi_3]\wedge d[\eta_1]
 -  \xi_1 \eta_3  d[\eta_1]\wedge d[\eta_2]
 +  \xi_2 \eta_3  d[\eta_1]\wedge d[\eta_1]$ & $(2, 1,-1)$ & $(45|0)$ & 1/0 \\
\hline
\end{tabular}}

\ssec{\protect$\fg=\fk\fa\fs(;1\xi)$} $2-d=0$. Nontrivial
cohomology are in degrees 0 and 1.

\nopagebreak

The weights are given relative to $\fo(6)$, but the corresponding vectors are
heighest only with respect to $\fo(4)=\fo(6)\cap\fg_0$.

\nopagebreak

Cohomology in $\deg= 0$,  $\dim=(6|6)$:

\vskip 0.2 cm

\nopagebreak

{\small

\renewcommand{\arraystretch}{1.4}\begin{tabular}{|c|c|c|c|c|}
\hline
N& $\fo(4)$-highest vectors & $\fo(6)$-weight & dim & mult \\
\hline
$[1]$ & $d[\xi_1 \xi_3]\wedge d[{{\xi }_1} {{\eta }_2}]$ & $(-2, 1, -1)$ & $(3|0)$ & 1/0 \\
\hline
$[2]$ & $d[{{\xi }_1} {{\eta }_2}]\wedge d[{{\xi }_1} {{\eta }_3}]$ & $(-2, 1, 1)$ &  $(3|0)$ & 1/0 \\
\hline
$[3]$ & ${{\xi }_1}d[{{\xi }_1} {{\xi }_3}]\wedge
d[{{\xi }_1} {{\eta }_2}]$ & $(-1, 1, -1)$ & $(0|3)$ & 2/1 \\
\hline
$[4]$ & ${{\xi }_1} d[{{\xi }_1} {{\eta }_2}]\wedge
d[{{\xi }_1} {{\eta }_3}]$ & $(-1, 1, -1)$ & $(0|3)$ & 2/1 \\
\hline
\end{tabular}}

\vskip 0.2 cm

$[2]\oplus [4]$ is $\fg_0$-irreducible;
[1] and [3] are glued as $\fg_0$-modules.

\vskip 0.2 cm

Cohomology in $\deg= 1$, $\dim=(8|8)$:

\vskip 0.2 cm

\nopagebreak

{\small
\renewcommand{\arraystretch}{1.4}\begin{tabular}{|c|c|c|c|c|}
\hline
N& $\fo(4)$-highest vectors & $\fo(6)$-weight & dim & mult \\
\hline
$[1]$ & $\xi_1 d[\xi_1 \eta_2] \wedge (\xi_2 d[\xi_3]-\eta_3 d[\eta_2])$ & $(0, 2, -1)$ & $(0|8)$ & 2/1 \\
\hline
$[2]$ & $\xi_1 d[\eta_2] \wedge (\xi_2 d[\xi_3]-\eta_3 d[\eta_2])$ & $(1, 2, -1)$ & $(8|0)$ & 1/0 \\
\hline
\end{tabular}}

\vskip 0.2 cm

A single irreducible $\fg_0$-module.

\ssec{\protect$\fg=\fk\fa\fs(; 3\eta)$} Cohomology in degrees 1 and 2.

$\deg =1$: a single irreducible $\fg_0$-module of $\dim=(12|12)$
(the vectors are highest with respect to
$\fgl(3)=\fo(6)\cap\fg_0$, the weights are the same as in
$\fo(6)$):

\vskip 0.2 cm

\nopagebreak

{\small
\renewcommand{\arraystretch}{1.4}\begin{tabular}{|c|c|c|c|c|}
\hline
N& $\fgl(3)$-highest vectors & weight & dim & mult \\
\hline
$[1]$ & $d[{{\eta }_1}{{\eta }_2}]\wedge d[{{\eta }_1}{{\eta }_3}]$ & $(2, 1, 1)$ & $(3|0)$ & 1/0 \\
\hline
$[2]$ & $d[{{\eta }_1}]\wedge d[{{\eta }_2}{{\eta }_3}]-d[{{\eta
}_2}]\wedge d[{{\eta }_1}{{\eta }_3}]+d[{{\eta }_3}]\wedge d[{{\eta
}_1}{{\eta }_2}]$ & $(1, 1, 1)$ & $(0|1)$ & 2/1 \\
\hline
$[3]$ & $d[{{\eta }_1}]\wedge d[{{\eta }_1}{{\eta }_2}]$ & $(2, 1, 0)$ & $(0|8)$ & 2/1 \\
\hline
$[4]$ & $d[1]\wedge d[{{\eta }_1}{{\eta }_2}]$ & $(1, 1, 0)$ & $(3|0)$ & 3/2 \\
\hline
$[5]$ & ${{d[{{\eta }_1}]}^2}$ & $(2, 0, 0)$ & $(6|0)$ & 3/2 \\
\hline
$[6]$ & $d[1]\wedge d[{{\eta }_1}]$ & $(1, 0, 0)$ & $(0|3)$ & 5/4 \\
\hline
\end{tabular}}

\vskip 0.2 cm

$\deg =2$: $\dim=(15|16)$, a single irreducible $\fg_0$-module:

\nopagebreak

{\small
\renewcommand{\arraystretch}{1.4}\begin{tabular}{|c|c|c|c|c|}
\hline
N& $\fgl(3)$-highest vectors & weight & dim & mult \\
\hline
$[1]$ & ${{\xi }_1}{{\eta }_3}d[{{\eta }_2}]\wedge
d[{{\eta }_1}{{\eta }_2}]-{{\xi }_2}{{\eta }_3}d[{{\eta }_1}]\wedge d[{{\eta }_1}{{\eta
}_2}]$ & $(2, 2, -1)$ & $(0|10)$ & 1/0 \\
\hline
$[2]$ & $\begin{matrix}-{{\xi }_1}{{\eta }_3}d[{{\eta }_1}]\wedge d[{{\eta }_2}]
+{{\xi }_2}{{\eta }_3}{{d[{{\eta }_1}]}^2}
+{{\xi }_1}{{\eta }_1}{{\eta }_3}d[{{\eta }_1}]\wedge d[{{\eta }_1}{{\eta }_2}]\\
+2{{\xi }_1}{{\eta }_2}{{\eta }_3}d[{{\eta }_2}]\wedge d[{{\eta }_1}{{\eta }_2}]
-{{\xi }_2}{{\eta }_2}{{\eta }_3}d[{{\eta }_1}]\wedge d[{{\eta }_1}{{\eta }_2}]\end{matrix}$
& $(2, 1, -1)$ & $(15|0)$ & 2/1 \\
\hline
$[3]$ & $\eta_3(\eta_1 d[\eta_1] + \eta_2 d[\eta_2])\wedge(\xi_1 d[\eta_2] - \xi_2 d[\eta_1])$ &  $(1, 1,-1)$ & $(0|6)$ & 3/2 \\
\hline
\end{tabular}}

\ssec{\protect$\fg=\fk\fa\fs(; 3\xi)$}
Cohomology in $\deg= 1$, $\dim=(52|52)$ (the vectors are highest with respect to
$\fgl(3)=\fo(6)\cap\fg_0$, the weights are given in the following table with respect to $\fo(6)$), where

$\fg_0$-modules:

$[A]=[2']+[3']+[4]+[5']+[6']+[7']$ of
$\dim=(12|12)$,  where $[2']$,  $[3']$,  $[5']$,  $[6']$,  $[7']$ are generated
by $[2_1]-[2_2]$, $[3_1]-2[3_2]+[3_3]$, $2[5_1]-[5_2]+[5_3]$, $[6_1]+[6_2]$, and
$[7_1]+[7_2]$, respectively;

$[B]=[8']+[9']+[10]+[11]+[12]+[13]$ of  $\dim=(24|24)$,
where $[8']$ and $[9']$ are generated by $[8_1]+2[8_2]$, and
$[9_1]+2[9_2]$, respectively;

$[C]=[B]+[3'']+[5'']+[6_1]+[7_1]+[8_1]+[9_1]$ of $\dim=(36|36)$,
where $[3'']$ and
$[5'']$ are generated by $[3_1]-[3_2]$ and $[5_2]-3[5_1]$, respectively;

$[D]=[A]+[C]$ of  $\dim=(48|48)$;

$[E]=[D]+[1]+[2'']+[3''']+[5''']$ of  $\dim=(52|52)$.

The modules $[A]$,  $[B]$,  $[C]/[B]$ and $[E]/[D]$ are irreducible,
$\dim([C]/[B])=(12|12)$ and  $\dim([E]/[D])=(4|4)$.

\vskip 0.2 cm

\nopagebreak

{\footnotesize
\renewcommand{\arraystretch}{1.4}\begin{tabular}{|c|c|c|c|c|}
\hline
N& $\fgl(3)$-highest vectors & weight & dim & mult \\
\hline
$[1]$ & ${{d[{{\xi }_1}{{\xi }_2}{{\xi }_3}]}^2}$ & $(-2, -2, -2)$ & $(1|0)$ &  1/0 \\
\hline
$[2_1]$ & $d[{{\xi }_2}{{\xi }_3}]\wedge d[{{\xi }_1}{{\xi }_2}{{\xi }_3}]$ & $(-1, -2, -2)$ & $2\cdot (0|3)$ & 2/0 \\
\cline{1-2}
$[2_2]$ & ${{\xi }_1} {{d[{{\xi }_1}{{\xi }_2}{{\xi }_3}]}^2}$ &&& \\
\hline
$[3_1]$ & $d[{{\xi }_3}]\wedge d[{{\xi }_1}{{\xi }_2}{{\xi }_3}]$ &&& \\
\cline{1-2}
$[3_2]$ & $d[{{\xi }_1}{{\xi }_3}]\wedge d[{{\xi }_2}{{\xi }_3}]$ & $(-1, -1, -2)$ & $3\cdot (3|0)$ & 4/1 \\
\cline{1-2}
$[3_3]$ & ${{\xi }_1}{{\xi }_2} {{d[{{\xi }_1}{{\xi }_2}{{\xi}_3}]}^2}$ &&& \\
\hline
$[4]$ & ${{\xi }_1} d[{{\xi }_2}{{\xi }_3}]\wedge
d[{{\xi }_1}{{\xi }_2}{{\xi }_3}]$ & $(0, -2, -2)$ & $(6|0)$ & 1/0 \\
\hline
$[5_1]$ & $d[1]\wedge d[{{\xi }_1}{{\xi }_2}{{\xi }_3}]$ &&& \\
\cline{1-2}
$[5_2]$ & $d[{{\xi }_1}]\wedge d[{{\xi }_2}{{\xi }_3}]-d[{{\xi
}_2}]\wedge d[{{\xi }_1}{{\xi }_3}]+d[{{\xi }_3}]\wedge d[{{\xi }_1}{{\xi }_2}]$ &$(-1, -1, -1)$ & $3\cdot(0|1)$ & 6/3 \\
\cline{1-2}
$[5_3]$ & ${{\xi }_1}{{\xi }_2}{{\xi }_3} {{d[{{\xi}_1}{{\xi }_2}{{\xi }_3}]}^2}$ &&& \\
\hline
$[6_1]$ & $d[{{\xi }_3}]\wedge d[{{\xi }_2}{{\xi }_3}]$ & $(0, -1, -2)$ & $2\cdot(0|8)$ & 4/2 \\
\cline{1-2}
$[6_2]$ & ${{\xi }_1}{{\xi }_2}d[{{\xi }_2}{{\xi }_3}]\wedge
d[{{\xi }_1}{{\xi }_2}{{\xi }_3}]$ &&& \\
\hline
$[7_1]$ & $d[1]\wedge d[{{\xi }_2}{{\xi }_3}]$ & $(0, -1, -1)$ & $2\cdot(3|0)$ & 7/5 \\
\cline{1-2}
$[7_2]$ & ${{\xi }_1}d[1]\wedge d[{{\xi }_1}{{\xi }_2}{{\xi }_3}]$ &&& \\
\hline
$[8_1]$ & ${{d[{{\xi }_3}]}^2}$ & $(0, 0, -2)$ & $2\cdot(6|0)$ & 4/2 \\
\cline{1-2}
$[8_2]$ & ${{\xi }_1}{{\xi }_2}d[{{\xi }_1}{{\xi }_3}]\wedge d[{{\xi }_2}{{\xi }_3}]$ &&& \\
\hline
$[9_1]$ & $d[1]\wedge d[{{\xi }_3}]$ & $(0, 0, -1)$ & $2\cdot(0|3)$ & 8/6 \\
\cline{1-2}
$[9_2]$ & ${{\xi }_1}{{\xi }_2}{{\xi }_3}d[{{\xi }_1}{{\xi }_3}]\wedge d[{{\xi
}_2}{{\xi }_3}]$ &&& \\
\hline
$[10]$ & ${{\xi }_1}{{d[{{\xi }_3}]}^2}$ & $(1, 0, -2)$ & $(0|15)$ & 2/1 \\
\hline
$[11]$ & ${{\xi }_1}d[1]\wedge d[{{\xi }_3}]$ & $(1, 0, -1)$ & $(8|0)$ & 4/3 \\
\hline
$[12]$ & ${{\xi }_1}{{\xi }_2}{{d[{{\xi }_3}]}^2}$ & $(1, 1, -2)$ & $(10|0)$ & 1/0 \\
\hline
$[13]$ & ${{\xi }_1}{{\xi }_2}d[1]\wedge d[{{\xi }_3}]$ & $(1, 1, -1)$ & $(0|6)$ & 1/0 \\
\hline
\end{tabular}\normalsize}

\vskip 0.2 cm

\ssec{\protect$\fg=\fm\fb(4|5)$} $2-d=0$. Here
$$
\text{$w(u_0)=0$, $w(u_1)=(1,0)$, $w(u_2)=(-1,1)$,
$w(u_3)=(0,-1)$, $w(\xi_i)=-w(u_i)$.}
$$

Cohomology: in $\deg=1$, hence, torsion free. A single irreducible
$\fg_0$-module of $\dim=(12|12)$ glued of the following
$\fsl(3)$-modules:

\vskip 0.2 cm

\nopagebreak

{\small
\renewcommand{\arraystretch}{1.4}\begin{tabular}{|c|c|c|c|c|}
\hline
N& $\fgl(3)$-highest vectors & weight & dim & mult \\
\hline
$[1]$ & ${u_0}d[{u_0}]\wedge d[{{\xi }_0}]$ & $(0, 0)$ & $(0|1)$ &  4/3 \\
\hline
$[2]$ & ${u_0}d[{u_0}]\wedge d[{u_3}]-{{\xi }_0}d[{u_3}]\wedge d[{{\xi }_0}]$ & $(0, 1)$ & $(3|0)$ & 6/5 \\
\hline
$[3]$ & ${u_0}d[{u_3}]\wedge d[{{\xi }_0}]$ & $(0, 1)$ & $(0|3)$ &  3/2 \\
\hline
$[4]$ & ${u_0}d[{{\xi }_0}]\wedge d[{{\xi }_1}]$ & $(1, 0)$ & $(3|0)$ & 3/2 \\
\hline
$[5]$ & ${u_0}d[{u_3}]\wedge d[{{\xi }_1}]+{u_1}d[{u_3}]\wedge d[{{\xi }_0}]$ & $(1, 1)$ & $(0|8)$ & 3/2 \\
\hline
$[6]$ & ${u_0}{{d[{u_3}]}^2}-{{\xi }_3}d[{u_3}]\wedge d[{{\xi }_0}]$ & $(0, 2)$ & $(6|0)$ & 1/0 \\
\hline
\end{tabular}}

\ssec{\protect$\fg=\fm\fb(4|5; 1)$} $2-d=0$.

Cohomology in $\deg =0$, $\dim=(22|22)$:

\vskip 0.2 cm

\nopagebreak

{\small
\renewcommand{\arraystretch}{1.4}\begin{tabular}{|c|c|c|c|c|}
\hline
N& $\fsl(2)\oplus\fsl(2)$-highest vectors & weight & dim & mult \\
\hline
$[1]$ & ${{\xi }_1}d[{u_2}{{\xi }_1}]\wedge d[{u_3}{{\xi }_1}]$ &  $(0, 0)$ & $(0|1)$ &  6/5 \\
\hline
$[2_1]$ & $d[{u_2}]\wedge d[{u_0}{u_3}-{{\xi }_1}{{\xi }_2}]-d[{u_3}]\wedge
d[{u_0}{u_2}+{{\xi }_1}{{\xi }_3}]$ & $(0, 1)$ & $2\times(2|0)$ & 6/4 \\
\cline{1-2}
$[2_2]$ & $d[{u_2}{{\xi }_1}]\wedge d[{u_3}{{\xi }_1}]$ &&& \\
\hline
$[3]$ & $d[{u_0}{u_2}+{{\xi }_1}{{\xi }_3}]\wedge
d[{u_3}{{\xi}_1}]-d[{u_0}{u_3}-{{\xi }_1}{{\xi }_2}]\wedge
d[{u_2}{{\xi }_1}]$ & $(0, 2)$ & $(0|3)$ & 2/1 \\
\hline
$[4]$ & $d[{{\xi }_2}]\wedge d[{u_3}{{\xi }_1}]$ & $(2, 1)$ & $(6|0)$ & 4/3 \\
\hline
$[5_1]$ & $d[{{\xi }_2}]\wedge d[{u_0}{u_3}-{{\xi }_1}{{\xi }_2}]$ & $(2, 2)$ & $2\times(0|9)$ & 3/1 \\
\cline{1-2}
$[5_2]$ & $d[{u_0}{u_3}-{{\xi }_1}{{\xi }_2}]\wedge d[{u_3}{{\xi }_1}]$ &&& \\
\hline
$[6]$ & ${{d[{u_0}{u_3}-{{\xi }_1}{{\xi }_2}]}^2}$ & $(2, 3)$ & $(12|0)$ & 1/0 \\
\hline
\end{tabular}}

\vskip 0.2 cm

$\fg_0$-modules:

$[A]=[4]+[5_1]+[5_2]+[6]$ of $\dim=18|18$;

$[B]=[2_2]+[3]+[A]$ of  $\dim=20|21$;

$[C]=[1]+[2_1]+[B]$ of  $\dim=22|22$.

Irreducible modules: $[A]$,  $[B]/[A]$ (of $\dim=2|3$) and $[C]/[B]$ (of $\dim=2|1$).

\ssec{$\fg=\fm\fb(4|5; K)$} $2-d=-1$.

The $\fg_0$-highest weight vectors are as follows:
$\fg_0=\fsl(2) \oplus \fgl(3)$; the first coordinate of the weight is given  with respect to $\fsl(2)$  realized as
$$
x_+ = q_0(\tau+q_0\xi_0-\sum_{i=1}^3q_i\xi_i)+2\xi_1\xi_2\xi_3,
$$
$x_-=\xi_0$; $\fgl(3)$ is realized as $x^i_j=q_i\xi_j$ ($i\neq j$)
and $x^i_i=\tau +q_i\xi_i-q_0\xi_0$.

In this realization
$$
\renewcommand{\arraystretch}{1.4}
\begin{array}{l}
\text{$w(q_0)=(1,-1,-1,-1)$, $w(u_1)=(-1,1,0,0)$,
$w(u_2)=(-1,0,1,0)$, $w(u_3)=(-1,0,0,1)$,}\\
\text{$w(\xi_0)=(-2,1,1,1)$, $w(\xi_1)=(0,-1,0,0)$,
$w(\xi_2)=(0,0,-1,0)$, $w(\xi_3)=(0,0,0,-1)$).}
\end{array}
$$

Cohomology: a single irreducible module in $\deg =-1$: \vskip 0.2
cm

{\small
\renewcommand{\arraystretch}{1.4}\begin{tabular}{|c|c|c|c|}
\hline
$\fsl(2)\oplus\fgl(3)$-highest vectors & weight & dim & mult \\
\hline
${u_0}{{dq}_3}\wedge {{dq}_3}$ & $(3, 0, 0, -2)$ & $(24|0)$ & 1/0 \\
\hline
\end{tabular}}

\ssec{$\fg=\fk\fsle(9|6; K)$} $2-d=0$. Cohomology are in degrees 0
and 1 and constitute irreducible modules. (The weights are given
in $\fgl$-basis of matrix diagonal units.)

\vskip 0.2 cm

\nopagebreak

{\small
\renewcommand{\arraystretch}{1.4}\begin{tabular}{|c|c|c|c|c|}
\hline
deg & $\fgl(5)$-highest vectors & $\fgl(5)$-weight & dim & mult \\
\hline
0 & ${\partial_5} {{d[\pi  {{dx}_4}{{dx}_5}]}^2}$ & $(0, 0, 0, -2, -3)$ &
$(175|0)$ & 1/0 \\
\hline 1 & $\begin{matrix} \sum {\partial_i}d[{{\pi
dx}_4}{{dx}_5}]\wedge d[{\partial_i}]+\\
 \sum{{\pi dx}_i}{{dx}_j}(d[{{\pi dx}_i}{{dx}_4}]\wedge d[{{\pi dx}_j}
{{dx}_5}]-\\
d[{{\pi dx}_i}{{dx}_5}]\wedge d[{{\pi dx}_j}{{dx}_4}]-\\ d[{{\pi
dx}_i}{{dx}_j}]\wedge d[{{\pi dx}_4}{{dx}_5}])
\end{matrix}$
& $(0, 0, 0, -1, -1)$ & $(0|10)$ & 3/2 \\
\hline
\end{tabular}

\ssec{$\fg=\fk\fs\fle(9|6; CK)$} $2-d=-1$. Cohomology are in degrees $-1$,  0, and 1.

Cohomology in $\deg =-1$,  $\dim=(36|36)$:

\vskip 0.2 cm

\nopagebreak

{\small
\renewcommand{\arraystretch}{1.4}\begin{tabular}{|c|c|c|c|c|}
\hline
N & $\fsl(2)\oplus\fsl(3)$-highest vectors & weight & dim & mult \\
\hline
$[1]$ & ${\partial_2}{{d[{{\pi dx}_2}{{dx}_5}]}^2}$ & $(3, 0, 2)$ &
$(24|0)$ & 1/0 \\
\hline
$[2]$ & ${\partial_2}d[{{\pi dx}_2}{{dx}_5}]\wedge d[{x_5}{\partial_1}]$ &
$(3, 0, 2)$ & $(0|24)$ & 1/0 \\
\hline
$[3]$ & ${\partial_2}d[{x_4}{\partial_1}]\wedge d[{x_5}{\partial_1}]$ &
$(3, 1, 0)$ & $(12|0)$ & 1/0 \\
\hline
$[4]$ & ${\partial_2}d[{{\pi dx}_2}{{dx}_4}]\wedge d[{x_5}
{\partial_1}]-{\partial_2}d[{{\pi dx}_2}{{dx}_5}]\wedge
d[{x_4}{\partial_1}]$ & $(3, 1, 0)$ & $(0|12)$ & 1/0 \\
\hline
\end{tabular}}

\vskip 0.2 cm

$\fg_0$-modules:

$[A]=[1]+[2]$ of  $\dim=24|24$;

$[B]=[A]+[3]+[4]$ of  $\dim=36|36$.

Irreducible modules: [A] and [B]/[A] (of $\dim=12|12$).

\vskip 0.2 cm

Cohomology in $\deg =0$, single $\fg_0$-module,  $\dim=10|10$:

\vskip 0.2 cm

\nopagebreak

{\small
\renewcommand{\arraystretch}{1.4}\begin{tabular}{|c|c|c|c|c|}
\hline
N & $\fsl(2)\oplus\fsl(3)$-highest vectors & weight & dim & mult \\
\hline
$[1]$ & ${\partial_5}d[{x_5}{\partial_1}]\wedge d[{x_5}
{\partial_2}]$ & $(0, 0, 3)$ & $(10|0)$ & 1/0 \\
\hline
$[2]$ & $\begin{matrix}{\partial_5}d[{{\pi dx}_1}
{{dx}_5}]\wedge d[{x_5}{\partial_1}]+{\partial_5}d[{{\pi dx}_2}
{{dx}_5}]\wedge d[{x_5}{\partial_2}]-\\
{{\pi dx}_3}{{dx}_4}d[{x_5}{\partial_1}]\wedge d[{x_5}{\partial_2}]
\end{matrix}$ & $(0, 0, 3)$ & $(0|10)$ & 1/0 \\
\hline
\end{tabular}}

\vskip 0.2 cm

Cohomology:  in $\deg =1$, a single $\fg_0$-module,  $\dim=(6|6)$:

\vskip 0.2 cm

\nopagebreak

{\small
\renewcommand{\arraystretch}{1.4}\begin{tabular}{|c|c|c|c|c|}
\hline
N & $\fsl(2)\oplus\fsl(3)$-highest vectors & weight & dim & mult \\
\hline
$[1]$ & $\begin{matrix} -{\partial_2}d[{x_5}{\partial_2}]
\wedge d[{\partial_1}]+\sum{\partial_i}d[{x_5}{\partial_1}]
\wedge d[{\partial_i}]+  \\
\sum\big({{\pi dx}_1}{{dx}_j}d\big[{{\pi dx}_i}{{dx}_j}\big]
\wedge d[{x_5}{\partial_i}]+
{x_j}{\partial_i}d\big[{x_j}{\partial_1}\big]\wedge d[{x_5}
{\partial_i}]\big)\end{matrix}$ & $(1,0, 1)$ & $(6|0)$ & 8/7 \\
\hline
$[2]$ & $\begin{matrix}-\sum{\partial_i}d[{{\pi dx}_2}
{{dx}_5}]\wedge d[{\partial_i}]+ \\
\sum{{\pi dx}_i}{{dx}_j}d[{{\pi dx}_2}{{dx}_5}]\wedge d\big[{{\pi
dx}_i}{{dx}_j}\big]- \\
{{\sum }_{1\leq i\leq 2}}{{\sum }_{3\leq j\leq
5}}{x_j}{\partial_i}d[{{\pi dx}_2}{{dx}_5}]\wedge d
\big[{x_j}{\partial_i}\big]\end{matrix}$
& $(1, 0, 1)$ & $(0|6)$ & 8/7 \\
\hline
\end{tabular}}
\vskip 0.2 cm

\nopagebreak

\ssec{$\fg=\fk\fsle(9|6)$} $2-d=0$. Cohomology in $\deg =1$, hence, torsion-free. $\dim=(168|167)$

\nopagebreak\vskip 0.2 cm\nopagebreak

{\tiny
\renewcommand{\arraystretch}{1.4}\begin{tabular}{|c|c|c|c|c|}
\hline
N & $\fgl(3)$-highest vectors & weight & dim & mult \\
\hline
$[1]$ &$\begin{matrix}\sum {{(-1)}^{p(\sigma )}}{\partial_{\sigma (1)}}
d[{{\pi dx}_{\sigma (2)}}{{dx}_{\sigma (3)}}]\wedge d[{x_{\sigma (4)}}{\partial_5}]- \\
2\sum {{\pi dx}_i}{{dx}_j}d[{x_i}{\partial_5}]\wedge d\big[{x_j}{\partial_5}\big]\end{matrix}$ & $(0,0,0)$ & $(0|1)$ & 2/1 \\
\hline
$[2_1]$ & $\sum {\partial_i}d[{\partial_i}]\wedge d[{x_4}{\partial_5}]+{\partial_4}d[{\partial_i}]\wedge
d[{x_i}{\partial_5}]$ & $(0,0,1)$ & $2\times(4|0)$ & 5/3 \\
\cline{1-2}
$[2_2]$ & $\begin{matrix}2{\partial_4}(-d[{{\pi dx}_1}{{dx}_2}]\wedge d[{{\pi dx}_3}{{dx}_4}]+
d[{{\pi dx}_1}{{dx}_3}]\wedge d[{{\pi dx}_2}{{dx}_4}]- \\
d[{{\pi dx}_1}{{dx}_4}]\wedge d[{{\pi dx}_2}{{dx}_3}])+
\sum {{\pi dx}_i}{{dx}_j}d\big[{{\pi dx}_i}{{dx}_j}\big]\wedge d[{x_4}{\partial_5}]
\end{matrix}$ &&&\\
\hline
$[3]$ & $\sum {\partial_4}d[{\partial_i}]\wedge d[{{\pi dx}_i}{{dx}_4}]- \sum
{{\pi dx}_{\sigma (1)}}{{dx}_{\sigma (2)}}d[{\partial_{\sigma (3)}}]\wedge
d[{x_4}{\partial_5}]$ & $(0,0,2)$ & $(0|10)$ & 2/1 \\
\hline
$[4_1]$ & ${{\sum }_{1\leq i\leq 4}}({\partial_i}d[{\partial_i}]\wedge d[{{\pi dx}_3}{{dx}_4}]+
{{\pi dx}_1}{{dx}_2}d[{\partial_i}]\wedge d[{x_i}{\partial_5}])$ & $(0,1,0)$ & $2\times(0|6)$ & 6/4 \\
\cline{1-2}
$[4_2]$ & $\begin{matrix}
\sum {{\pi dx}_i}{{dx}_i}d\big[{{\pi dx}_i}{{dx}_j}\big]\wedge d[{{\pi dx}_3}{{dx}_4}]-\\
2\sum{{(-1)}^{\sign(\sigma )}}{{\pi dx}_1}{{dx}_2} +
d[{{\pi dx}_1}{{dx}_{\sigma (2)}}]\wedge d[{{\pi dx}_{\sigma (3)\sigma (4)}}]
\end{matrix}$ &&&\\
\hline
$[5_1]$ & ${\partial_4}d[{x_3}{\partial_5}]\wedge d[{x_4}{\partial_5}]$ & $(0, 1, 1)$ & $2\times(20|0)$& 4/2 \\
\cline{1-2}
$[5_2]$ & $-2{\partial_4}d[{\partial_1}]\wedge d[{\partial_2}]- {x_1}{\partial_5}d[{\partial_2}]\wedge d[{x_4}{\partial_5}]+
{x_2}{\partial_5}d[{\partial_1}]\wedge d[{x_4}{\partial_5}]$ &&& \\
\hline
$[6]$ & ${\partial_4}d[{{\pi dx}_3}{{dx}_4}]\wedge d[{x_4}{\partial_5}]$ & $(0, 1, 2)$ & $(0|45)$ & 1/0 \\
\hline
$[7]$ & $\begin{matrix}{{\pi dx}_1}{{dx}_2}d[{x_3}{\partial_5}]\wedge d[{x_4}{\partial_5}]-
{\partial_3}d[{{\pi dx}_3}{{dx}_4}]\wedge d[{x_4}{\partial_5}]+ \\
{\partial_4}d[{{\pi dx}_3}{{dx}_4}]\wedge d[{x_3}{\partial_5}]\end{matrix}$ & $(0, 2, 0)$ & $(0|20)$ & 2/1 \\
\hline
$[8]$ & ${{\pi dx}_1}{{dx}_2}d[{{\pi dx}_3}{{dx}_4}]\wedge d[{x_4}{\partial_5}]-
{\partial_4}{{d[{{\pi dx}_3}{{dx}_4}]}^2}$ & $(0, 2, 1)$ & $(60|0)$ & 1/0 \\
\hline
$[9]$ & $\sum ({\partial_i}d[{\partial_1}]\wedge d[{\partial_i}]+{x_1}{\partial_5}d[{\partial_i}]\wedge
d[{x_i}{\partial_5}])$ & $(1, 0, 0)$ & $(4|0)$ & 5/4 \\
\hline
$[10]$ & $\begin{matrix}\sum {(-1)}^{p(\sigma )}({\partial_{\sigma (2)}}d[{{\pi dx}_{\sigma
(3)}} {{dx}_{\sigma (4)}}]\wedge d[{x_4}{\partial_5}]+ \\
{\partial_4}d[{{\pi dx}_{\sigma (2)}}{{dx}_{\sigma (3)}}]\wedge d[{x_{\sigma (3)}}{\partial_5}])\end{matrix}$ &
$(1, 0, 1)$ & $(0|15)$ & 4/3 \\
\hline
$[11]$ & ${\partial_4}d[{\partial_1}]\wedge d[{x_4}{\partial_5}]$ & $(1, 0, 2)$ & $(36|0)$ & 4/3 \\
\hline
$[12]$ & $2{x_1}{\partial_5}d[{x_3}{\partial_5}]\wedge d[{x_4}{\partial_5}]-{\partial_3}d[{\partial_1}]\wedge d[{x_4}{\partial_5}]+
{\partial_4}d[{\partial_1}]\wedge d[{x_3}{\partial_5}]$ & $(1, 1, 0)$ & $(20|0)$ & 4/3 \\
\hline
$[13]$ & ${{\pi dx}_1}{{dx}_2}d[{\partial_1}]\wedge d[{x_4}{\partial_5}]+{\partial_4}d[{\partial_1}]\wedge
d[{{\pi dx}_3}{{dx}_4}]$ & $(1, 1, 1)$ & $(0|64)$ & 3/2 \\
\hline
\end{tabular}}

\vskip 0.2cm

$\fg_0$-modules:

$[A]=[2']+[4']+[9]$,  $\dim=(8|6)$,  where $[2']$ and $[4']$ are generated by $2[2_1]-[2_2]$ and
$2[4_1]-[4_2]$, respectively;

[B]=[A]+[1],  $\dim=(8|7)$;

$[C]=[1]+\dots+[13]$,  $\dim=(168|167)$.

The modules [A],  [B]/[A] (of $\dim=(0|1)$), and [C]/[B]
(of $\dim=(160|160)$) are irreducible.

\ssec{$\fg=\fk\fs\fle(9|6; 2)$} $2-d=0$. Cohomology are in degrees 0 and 1.

Cohomology in $\deg =0$, $\dim=(140|140)$:

\nopagebreak\vskip 0.2 cm\nopagebreak

{\tiny
\renewcommand{\arraystretch}{1.4}\begin{tabular}{|c|c|c|c|c|}
\hline
N & $\fsl(3)\oplus\fsl(2)$-highest vectors & weight & dim & mult \\
\hline
$[1]$ & ${{\pi dx}_4}{{dx}_5}{{d[{{\pi x}_5}{{dx}_4}{{dx}_5}]}^2}$ & $(0,0,2)$ & $(0|3)$ & 2/1 \\
\hline
$[2]$ & ${\partial_3}d[{\partial_4}]\wedge d[{{\pi x}_4}{{dx}_4}{{dx}_5}]+{\partial_3}d[{\partial_5}]\wedge
d[{{\pi x}_5}{{dx}_4}{{dx}_5}]$
& $(0,1,0)$ & $(0|3)$ & 4/3 \\
\hline
$[3_1]$ & ${\partial_3}d[{\partial_4}]\wedge d[{{\pi x}_5}{{dx}_4}{{dx}_5}]$ & $(0,1,2)$ & $2\cdot(0|9)$ & 5/3 \\
\cline{1-2}
$[3_2]$ & ${{\pi dx}_4}{{dx}_5}d[{x_5}{\partial_1}]\wedge d[{x_5}{\partial_2}]$ & & & \\
\hline
$[3_3]$ & ${\partial_3}{{d[{{\pi x}_5}{{dx}_4}{{dx}_5}]}^2}$ & $(0,1,2)$ & $(9|0)$ & 2/1 \\
\hline
$[4]$ & $\begin{matrix}{\partial_3}d[{{\pi dx}_3}{{dx}_4}]\wedge d[{{\pi x}_5}{{dx}_4}{{dx}_5}] - \\
{\partial_3}d[{{\pi dx}_3}{{dx}_5}]\wedge d[{{\pi x}_4}{{dx}_4}{{dx}_5}]\end{matrix}$ & $(0, 2, 0)$ & $(6|0)$ & 1/0 \\
\hline
$[5_1]$ & ${\partial_3}d[{\partial_4}]\wedge d[{{\pi dx}_3}{{dx}_5}]$
& $(0, 2, 2)$ & $(0|18)$ & 2/0 \\
\hline 
$[5_2]$ & ${\partial_3}d[{{\pi dx}_3}{{dx}_5}]\wedge
d[{{\pi x}_5}{{dx}_4}{{dx}_5}]$ &$(0, 2, 2)$ & $2\cdot(18|0)$ & 2/1 \\
\cline{1-2}
$[5_3]$ & ${\partial_3}d[{x_5}{\partial_1}]\wedge d[{x_5}{\partial_2}]$ &&&\\
\hline
$[6]$ & ${\partial_3}{{d[{{\pi dx}_3}{{dx}_5}]}^2}$ &  $(0, 3, 2)$ & $(30|0)$ & 1/0 \\
\hline
$[7]$ & $\begin{matrix}{{\pi dx}_4}{{dx}_5}d[{x_4}{\partial_1}]\wedge d[{{\pi x}_5}{{dx}_4}{{dx}_5}]- \\
{{\pi dx}_4}{{dx}_5}d[{x_5}{\partial_1}]\wedge d[{{\pi x}_4}{{dx}_4}{{dx}_5}]\end{matrix}$ & $(1, 0, 0)$ & $(3|0)$ & 3/2 \\
\hline
$[8]$ & ${{\pi dx}_4}{{dx}_5}d[{x_5}{\partial_1}]\wedge d[{{\pi x}_5}{{dx}_4}{{dx}_5}]$
& $(1, 0, 2)$ & $(9|0)$ & 3/2 \\
\hline
$[9_1]$ & ${\partial_3}d[{\partial_4}]\wedge d[{x_4}{\partial_1}]+{\partial_3}d[{\partial_5}]\wedge d[{x_5}{\partial_1}]$ &
$(1, 1, 0)$ & $(8|0)$ & 3/2 \\
\hline
$[9_2]$ & ${\partial_3}d[{x_4}{\partial_1}]\wedge d[{{\pi x}_5}{{dx}_4}{{dx}_5}]-
{\partial_3}d[{x_5}{\partial_1}]\wedge d[{{\pi x}_4}{{dx}_4}{{dx}_5}]$ & $(1, 1, 0)$ & $(0|8)$ & 1/0 \\
\hline
$[10_1]$ & ${\partial_3}d[{\partial_4}]\wedge d[{x_5}{\partial_1}]$ &  $(1, 1, 2)$ & $(24|0)$ & 3/2 \\
\hline
$[10_2]$ & ${\partial_3}d[{x_5}{\partial_1}]\wedge d[{{\pi x}_5}{{dx}_4}{{dx}_5}]$ &  $(1, 1, 2)$ & $(0|24)$ & 1/0 \\
\hline
$[11]$ & ${\partial_3}d[{{\pi dx}_3}{{dx}_4}]\wedge d[{x_5}{\partial_1}]-
{\partial_3}d[{{\pi dx}_3}{{dx}_5}]\wedge d[{x_4}{\partial_1}]$ & $(1, 2, 0)$ & $(0|15)$ & 1/0 \\
\hline
$[12]$ & ${\partial_3}d[{{\pi dx}_3}{{dx}_5}]\wedge d[{x_5}{\partial_1}]$ & $(1, 2, 2)$ & $(0|45)$ & 1/0 \\
\hline
$[13]$ & ${{\pi dx}_4}{{dx}_5}d[{x_4}{\partial_1}]\wedge d[{x_5}{\partial_1}]$ & $(2, 0, 0)$ & $(0|6)$ & 2/1 \\
\hline
$[14]$ & ${\partial_3}d[{x_4}{\partial_1}]\wedge d[{x_5}{\partial_1}]$ & $(2, 1, 0)$ & $(15|0)$ & 1/0 \\
\hline
\end{tabular}}

\vskip 0.2 cm

{\footnotesize The column \lq\lq mult" shows the multiplicity of the homogeneous
(even or odd) component of given weight.}

\vskip 0.2 cm

$\fg_0$-modules:

$[A]=[3_1-3_2]+[5_1]+[5_2+5_3]+[6]+[10_1]+[12]$ of  $\dim=72|72$;

$[B]=[A]+[2]+[4]+[7]+[9_1]+[9_2]+[11]+[13]+[14]$ of  $\dim=104|104$;

$[C]=\text{all}$.

The modules [A],  [B]/[A] (of $\dim=(32|32)$), and [C]/[B]
(of $\dim=(36|36)$) are irreducible.

\vskip 0.2 cm

Cohomology  in $\deg =1$, $\dim=(8|8)$: the
$\fsl(3)\oplus\fsl(2)$-highest weights are as follows:
even ones are (0, 0, 1), (1, 0, 1); odd ones are
(0, 0, 1), (0, 1, 1). The corresponding highest weight vectors are
too complicated to be included here.

\Ack P.G., D.~L. and I.~Shch.  acknowledge financial support of TBSS, Stockholm;
Universit\'e Marseille-Aix  and MPIM, Bonn, where the final molding had been performed; and RFBR grant 01-01-00490a, respectively.

\end{document}